\documentclass[11pt]{amsart}

\usepackage[utf8]{inputenc}
\usepackage[T1]{fontenc}
\usepackage{amsmath,amssymb,amsthm}
\usepackage[margin=1.1in]{geometry}

\newtheorem{theorem}{Theorem}[section]
\newtheorem{lemma}[theorem]{Lemma}
\newtheorem{proposition}[theorem]{Proposition}
\newtheorem{corollary}[theorem]{Corollary}
\theoremstyle{definition}
\newtheorem{definition}[theorem]{Definition}
\newtheorem{example}[theorem]{Example}
\newtheorem{question}[theorem]{Question}
\theoremstyle{remark}
\newtheorem{remark}[theorem]{Remark}

\newcommand{\C}{\mathbb{C}}
\newcommand{\Pj}{\mathbb{P}}
\newcommand{\Hm}{\mathcal{H}}
\newcommand{\Vol}{\operatorname{Vol}}
\newcommand{\Sing}{\operatorname{Sing}}
\newcommand{\dc}{d^{c}}
\newcommand{\ddc}{dd^{c}}
\newcommand{\bd}{\bar\partial}
\newcommand{\ovl}[1]{\overline{#1}}
\newcommand{\nrm}[1]{\lVert #1\rVert}
\newcommand{\tab}[2]{t^{#1\bar{#2}}}
\newcommand{\dashto}{\dashrightarrow}

\numberwithin{equation}{section}
\begin{document}

\title[Extension of meromorphic mappings with finite volume]{On the Extension of Meromorphic Mappings and Correspondences with Finite Volume}

\author{Fethi Neji}
\address{University of Gabes, Higher Institute of Informatics of Medenine,
, LR17ES11 Mathematics and Applications
laboratory, 6072, Gabes,
}
\email{fethi.neji@isim.rnu.tn}
\date{October 23, 2026}

\begin{abstract}
Let $X$ be a compact complex manifold. Using Bishop's extension theorem, we show that a meromorphic map, or more generally a meromorphic correspondence, $f\colon\Delta^k_*\dashto X$ extends meromorphically across $0$ if and only if its graph has finite $2k$-dimensional volume near $\{0\}\times X$; equivalently, all mixed integrals $\int\omega_E^{k-m}\wedge(f^*\omega)^m$, $1\le m\le k$, are finite. The criterion extends to arbitrary proper analytic subsets, and yields an energy criterion with a sharp exponent, a dichotomy for maps into Hopf manifolds, a rationality criterion for maps between projective manifolds, and, for $k=2$ and pluriclosed targets, the equivalence of extension with finiteness of the single integral $\int(f^*\omega)^2$.
\end{abstract}


\maketitle

\section{Introduction}\label{sec:intro}

\subsection*{The extension problem}
A classical theorem of Hartogs asserts that, for $k\ge2$, every holomorphic function on the punctured polydisc $\Delta^k_*=\Delta^k\setminus\{0\}$ extends holomorphically to $\Delta^k$. We study the analogue for meromorphic maps $f\colon\Delta^k_*\dashto X$ with values in a complex manifold $X$, and more generally for meromorphic correspondences. For $k=1$ there is no extension phenomenon at all: $f(z)=[1/z]\in\C/\Lambda$ has an essential singularity at $0$ (Remark~\ref{rem:k1}). For $k\ge2$ and $X$ compact, the answer depends on the geometry of $X$: it is positive for K\"ahler targets \cite{Siu}, but extension may fail for non-K\"ahler targets, for instance Hopf surfaces (see \cite{Iv04} and Example~\ref{ex:hopf} below). It is therefore natural to look for a condition on the map $f$ itself that forces extension. The present work originates in an open question raised by Ivashkovich in his survey on extension properties of equidimensional holomorphic mappings \cite{Iv08}, which led us to consider the volume of the graph as the natural condition.

The condition studied here is the finiteness of the volume of the graph, counted with multiplicity. It is natural for the following reason. By a theorem of Bishop \cite{Bis,Chi}, the closure of an analytic set of locally finite volume is analytic, and conversely, by Lelong's theorem \cite{Lel}, every analytic set has locally finite volume. If $f$ extends meromorphically across $0$, its graph is therefore an analytic subset of $\Delta^k\times X$ and has finite volume near $\{0\}\times X$. The point of this paper is that, for a compact target, the converse holds.

\subsection*{Setting and definitions}
Throughout, $\Delta\subset\C$ is the unit disc, $\Delta^k$ the unit polydisc, $\omega_E=\frac i2\sum_j dz_j\wedge d\bar z_j$ the Euclidean K\"ahler form on $\C^k$ and $dV$ Lebesgue measure. A Hermitian metric on a complex manifold is identified with its fundamental positive $(1,1)$-form $\omega$. We use the following notions; the formal definitions are in Section~\ref{sec:prelim}.

\smallskip\noindent\emph{Meromorphic maps.} Let $D$ and $X$ be complex manifolds. A \emph{meromorphic map} $f\colon D\dashto X$ is given by an irreducible analytic subset $\Gamma_f\subset D\times X$, its graph, such that the projection $\pi\colon\Gamma_f\to D$ is proper and is a proper modification: $\pi$ is biholomorphic over $D\setminus S$, where $S\subset D$ is a nowhere dense analytic set with $\pi^{-1}(S)$ nowhere dense in $\Gamma_f$. The \emph{indeterminacy set} is $I_f=\{z\in D:\dim\pi^{-1}(z)\ge1\}$; it is analytic of codimension $\ge2$, and $f$ is holomorphic on $D\setminus I_f$ (Definition~\ref{def:mero}).

\smallskip\noindent\emph{Meromorphic correspondences.} Let $D$ be connected of dimension $k$. A \emph{meromorphic correspondence} $F\colon D\dashto X$ \emph{of degree $d$} is a pure $k$-dimensional reduced analytic subset $\Gamma_F\subset D\times X$ such that $\pi\colon\Gamma_F\to D$ is proper and is an unramified $d$-sheeted covering over $D\setminus S$, for a nowhere dense analytic set $S\subset D$ with $\pi^{-1}(S)$ nowhere dense in $\Gamma_F$ (Definition~\ref{def:corr}). A meromorphic map is a correspondence of degree $1$. Unlike a map, a correspondence may have a fibre of positive dimension over a point of $D$ (Example~\ref{ex:corr}).

\smallskip\noindent\emph{Volume of the graph and mixed integrals.} Let $X$ carry a Hermitian metric $\omega$ and give $\Delta^k\times X$ the product metric. For a pure $k$-dimensional analytic set $A\subset\Delta^k\times X$, Wirtinger's equality gives $\operatorname{Vol}(A)=\mathcal H^{2k}(A)=\int_{A_{\rm reg}}(\omega_E+\omega)^k/k!$, and Lelong's theorem says that this is finite on compact sets. For a meromorphic map $f$ on $\Delta^k_*$ and $U\Subset\Delta^k$, the volume of the graph expands as
\[
\operatorname{Vol}\bigl(\Gamma_f\cap(U\times X)\bigr)=\sum_{m=0}^{k}\frac1{m!\,(k-m)!}V_m(U),\qquad
V_m(U)=\int_{U\setminus(\{0\}\cup I_f)}\omega_E^{k-m}\wedge(f^*\omega)^m,
\]
with $V_0(U)=\int_U\omega_E^k<\infty$ (Lemma~\ref{lem:expansion}); here $V_m$ is the $m$-th \emph{mixed integral}, and $V_k=\int(f^*\omega)^k$ is the integral of the pull-back of the volume form of the image. If $X$ is compact, any two Hermitian metrics are comparable, so the finiteness of these quantities does not depend on $\omega$ (Lemma~\ref{lem:volume-general}).

\smallskip\noindent\emph{Cluster set.} The cluster set of $f$ at $0$ is $\{x\in X:(0,x)\in\overline{\Gamma_f}\}$, the closure being taken in $\Delta^k\times X$; it is the set of limits of $f(z_j)$ along sequences $z_j\to0$ in $\Delta^k_*\setminus I_f$. Boundedness of the cluster set, in the sense of Theorem~\ref{thm:bounded}, replaces compactness of $X$.

\smallskip\noindent\emph{Hermitian geometry.} Let $d^c=\frac i2(\bar\partial-\partial)$, so that $dd^c=i\partial\bar\partial$. A Hermitian metric $\omega$ is \emph{pluriclosed} (or SKT) if $dd^c\omega=0$. When $\dim X=2$ this is the Gauduchon condition, and every compact complex surface carries such a metric by \cite{Gau}.

\smallskip\noindent\emph{Hopf manifolds.} For $n\ge2$ and a linear contraction $\gamma$ of $\C^n$, the \emph{Hopf manifold} is $H=(\C^n\setminus\{0\})/\langle\gamma\rangle$; for $n=2$ it is a Hopf surface. It is compact, diffeomorphic to $S^1\times S^{2n-1}$, hence not K\"ahler, and the quotient map $p\colon\C^n\setminus\{0\}\to H$ is its universal covering.

\smallskip\noindent\emph{Energy.} For a holomorphic map $f$ between Hermitian manifolds, $\lVert df\rVert$ denotes the operator norm of the differential; the integral $\int\lVert df\rVert^{2n}dV$, $n=\dim$ of the source, is the energy criterion of Corollary~\ref{cor:energy}.

\subsection*{Main results}
Our first three results are formulated for the punctured polydisc.

\begin{theorem}[Maps]\label{thm:maps}
Let $X$ be a compact complex manifold with a Hermitian metric $\omega$, $k\ge1$, and $f\colon\Delta^k_*\dashto X$ meromorphic with graph $\Gamma_f$. Assume
\begin{equation}\label{eq:fv}
\Vol\bigl(\Gamma_f\cap(U\times X)\bigr)<\infty
\end{equation}
for some neighbourhood $U$ of $0$, where $\Vol$ is the $2k$-dimensional Hausdorff measure for the product metric. Then $f$ extends to a meromorphic map $\Delta^k\dashto X$. For $U\Subset\Delta^k$, condition \eqref{eq:fv} is equivalent to $V_m(U)<\infty$ for $1\le m\le k$, where
\[
V_m(U)=\int_{U\setminus(\{0\}\cup I_f)}\omega_E^{k-m}\wedge(f^*\omega)^m,\qquad \omega_E=\frac i2\sum_j dz_j\wedge d\bar z_j .
\]
\end{theorem}

\begin{theorem}[Correspondences]\label{thm:corr}
Let $X$ be compact, $k\ge1$, and $F\colon\Delta^k_*\dashto X$ a meromorphic correspondence of degree $d$ (Definition~\ref{def:corr}) whose graph satisfies \eqref{eq:fv}. Then $F$ extends to a meromorphic correspondence $\widehat F\colon\Delta^k\dashto X$ of degree $d$. The fibre of $\widehat F$ over $0$ may have positive dimension.
\end{theorem}

\begin{theorem}[Bounded cluster set]\label{thm:bounded}
Let $X$ be any complex manifold and $f$ (or $F$) as above satisfying \eqref{eq:fv}. Assume that the graph over $U\setminus\{0\}$ lies in $(U\setminus\{0\})\times K$ for some compact $K\subset X$. Then the conclusions of Theorems~\ref{thm:maps} and~\ref{thm:corr} hold.
\end{theorem}

The next two results are consequences and refinements of the first three: a characterization together with a criterion in terms of energy (Theorem~\ref{thm:char}), and a sharper statement in dimension two (Theorem~\ref{thm:pluri}).

\begin{theorem}[Characterization, general exceptional sets, energy]\label{thm:char}
Let $X$ be compact.
\begin{enumerate}
\item $f\colon\Delta^k_*\dashto X$ extends meromorphically if and only if \eqref{eq:fv} holds (Theorem~\ref{thm:maps}, Proposition~\ref{prop:converse}).
\item If $Z$ is a proper closed analytic subset of a connected complex manifold $Y$, a meromorphic map $Y\setminus Z\dashto X$ extends to $Y$ if and only if its graph has finite volume over every compact subset of $Y$ (Theorem~\ref{thm:general}).
\item If $Z\subset Y$ is a closed analytic subset of codimension $c\ge1$ and $\int_{K\setminus Z}\nrm{df}^{2n}dV<\infty$ for every compact $K\subset Y$ ($n=\dim Y$), then $f$ extends. If $Z$ is a complex submanifold, this holds when $\nrm{df}=O(\operatorname{dist}(\cdot,Z)^{-a})$ with $a<c/n$. For $Z=\{0\}\subset\Delta^n$ this reads $a<1$, and $a=1$ is attained by maps without extension (Corollary~\ref{cor:energy}, Remark~\ref{rem:sharp}).
\end{enumerate}
\end{theorem}

\begin{theorem}[Pluriclosed targets, $k=2$]\label{thm:pluri}
Let $X$ be a compact complex manifold with a Hermitian metric $\omega$ such that $\ddc\omega=0$ (this is the Gauduchon condition when $\dim X=2$, so by \cite{Gau} the hypothesis holds on every compact complex surface). Let $f\colon\Delta^2_*\to X$ be holomorphic. Then $V_1(U)<\infty$ for some neighbourhood $U$ of $0$, and
\[
f\ \text{extends meromorphically across }0\iff\int_{U\setminus\{0\}}(f^*\omega)^2<\infty .
\]
\end{theorem}

\subsection*{Comments on the results}
Theorem~\ref{thm:maps} is a characterization: by Proposition~\ref{prop:converse} the converse holds, so for a compact target a meromorphic map on $\Delta^k_*$ extends if and only if \eqref{eq:fv} holds. Condition \eqref{eq:fv} is not automatic, even for compact targets: for the Hopf surface of Example~\ref{ex:hopf} the integral $V_2$ diverges logarithmically and $f$ has no extension, whereas Proposition~\ref{prop:stein} shows that for $k\ge2$ no volume hypothesis is needed for holomorphic maps when the universal cover of $X$ is Stein. For correspondences the extended graph may contain a whole fibre over $0$ (Example~\ref{ex:corr}), so the extension need not be finite over $0$. The role of the mixed integrals is clarified by Remark~\ref{rem:Vk}: no inequality $V_1\le C\,V_k$ holds uniformly, and whether finiteness of $V_k$ alone forces extension is open in general (Question~\ref{q:single}); Theorem~\ref{thm:pluri} answers it when $k=2$, for holomorphic maps and pluriclosed targets.

\subsection*{Ideas of the proofs}
By Lemma~\ref{lem:expansion}, finiteness of the volume of the graph amounts to finiteness of $V_1,\dots,V_k$. Away from $\{0\}\times X$ the graph is analytic, so Lelong's theorem gives local finiteness there (Lemma~\ref{lem:localfin}); the hypothesis \eqref{eq:fv} supplies it near $\{0\}\times X$. Bishop's theorem then shows that the closure of the graph is analytic, pure $k$-dimensional and has no component over $\{0\}$ (Lemma~\ref{lem:closure}). Properness of the projection comes from the compactness of $X$ (or from the bounded cluster set), and Lemma~\ref{lem:tomero}, which uses the theorems of Remmert and Grauert, identifies the closure as the graph of a meromorphic map. For correspondences, one controls the ramification and the number of sheets (Section~\ref{sec:corr}). The energy criterion follows from the density bound of Lemma~\ref{lem:density}, and the statements on Hopf manifolds and on targets with Stein universal cover combine lifting to the universal cover with Hartogs' theorem. The proof of Theorem~\ref{thm:pluri} is different in nature: for a pluriclosed target the form $f^*\omega$ satisfies an identity (Lemma~\ref{lem:pluri}) which makes the slice area $a(z_1)=\int_{\{z_1\}\times\Delta(r_1)}f^*\omega$ solve a Poisson equation with smooth right-hand side; B\^ocher's theorem then gives $a(z_1)=O(\log\tfrac1{|z_1|})$, which is integrable, so $V_1<\infty$ (Section~\ref{sec:pluri}).

\subsection*{Relation to previous work}
The analytic input is classical: the closure of an analytic set of locally finite volume is analytic (Bishop \cite{Bis}, Chirka \cite{Chi}), and conversely analytic sets have locally finite volume (Lelong \cite{Lel}). Extension of meromorphic maps into K\"ahler manifolds is due to Siu \cite{Siu}; extension problems for non-K\"ahler targets, and the notion of disk-convexity, are studied in \cite{Iv04}; the question addressed here was motivated by an open problem from \cite{Iv08}. The contribution here is the systematic use of the equivalence above as an if-and-only-if extension criterion, with the consequences listed above (mixed-integral formulation, correspondences, energy criterion, Hopf dichotomy). The new analytic step is Theorem~\ref{thm:pluri}, which uses the pluriclosedness of $f^*\omega$ to control $V_1$ by a slice argument. In Section~\ref{sec:applications} we record further consequences: the criterion over an arbitrary base $Y\setminus Z$ (Theorem~\ref{thm:corr-general}), a rationality criterion for maps between projective manifolds (Corollary~\ref{cor:rational}), and a statement for families of correspondences (Proposition~\ref{prop:families}). Whether $V_k$ alone suffices in general is left open (Section~\ref{sec:open}).

\subsection*{Organization of the paper}
Section~\ref{sec:prelim} fixes the definitions and the volume computations. Section~\ref{sec:bishop} recalls Bishop's theorem. Section~\ref{sec:maps} proves Theorems~\ref{thm:maps} and~\ref{thm:bounded}, and Section~\ref{sec:corr} proves Theorem~\ref{thm:corr}. Section~\ref{sec:charac} contains the converse, the general exceptional sets, the energy criterion, the Hopf and Stein cases and an example of a correspondence with a positive-dimensional fibre. Section~\ref{sec:pluri} is devoted to pluriclosed targets in dimension two (Theorem~\ref{thm:pluri}). Section~\ref{sec:applications} gives further consequences and Section~\ref{sec:open} lists open questions.

\section{Preliminaries}\label{sec:prelim}

\begin{definition}[Meromorphic map]\label{def:mero}
Let $D$, $X$ be complex manifolds. A meromorphic map $f\colon D\dashto X$ is given by an irreducible analytic subset $\Gamma_f\subset D\times X$ such that the projection $\pi\colon\Gamma_f\to D$ is proper and a proper modification: there is a nowhere dense analytic set $S\subset D$ such that $\pi^{-1}(S)$ is nowhere dense in $\Gamma_f$ and $\pi$ is biholomorphic over $D\setminus S$. The indeterminacy set is $I_f=\{z:\dim\pi^{-1}(z)\ge1\}$; it is analytic of codimension $\ge2$ and $f$ is holomorphic on $D\setminus I_f$.
\end{definition}

\begin{definition}[Meromorphic correspondence]\label{def:corr}
Let $D$ be a connected complex manifold, $k=\dim D$. A meromorphic correspondence $F\colon D\dashto X$ is a pure $k$-dimensional reduced analytic subset $\Gamma_F\subset D\times X$ such that $\pi\colon\Gamma_F\to D$ is proper and there are a nowhere dense analytic set $S\subset D$ with $\pi^{-1}(S)$ nowhere dense in $\Gamma_F$ and an integer $d\ge1$ such that $\pi$ is an unramified $d$-sheeted covering over $D\setminus S$. By Lemma~\ref{lem:conn}, $D\setminus S$ is connected, so $d$ is independent of $S$.
\end{definition}

\subsection{Volume}
Let $\omega$ be a smooth Hermitian metric on $X$ and $\Omega=\omega_E+\omega$ on $\C^k\times X$ (or $\pi^*\omega_E+p^*\omega$ on $D\times X$). For a pure $k$-dimensional analytic set $A$, Wirtinger's equality is a pointwise algebraic identity for complex $k$-planes in a Hermitian space, so $\Hm^{2k}(A)=\int_{A_{\mathrm{reg}}}\Omega^k/k!$, since $A_{\mathrm{sing}}$ has dimension $<k$; see \cite{Chi,Dem,Fed}. By Lelong's theorem \cite{Lel}, a pure $k$-dimensional analytic subset of a complex manifold has locally finite $2k$-dimensional volume, i.e.\ finite volume on every compact set.

\begin{lemma}[Expansion of the graph volume]\label{lem:expansion}
For $f\colon\Delta^k_*\dashto X$ meromorphic and $U\Subset\Delta^k$ open,
\[
\Vol\bigl(\Gamma_f\cap(U\times X)\bigr)=\sum_{m=0}^{k}\frac1{m!\,(k-m)!}\,V_m(U),\qquad V_m(U)\in[0,\infty],
\]
with $V_0(U)=\int_U\omega_E^k<\infty$. Hence the volume is finite if and only if $V_1,\dots,V_k$ are finite.
\end{lemma}

\begin{proof}
Over $\Delta^k_*\setminus I_f$, $\Gamma_f$ is the graph of $f$, a complex submanifold; the part over $I_f$ is a proper analytic subset of the irreducible set $\Gamma_f$, hence of zero $\Hm^{2k}$-measure. Now $\Omega=\pi^*\omega_E+p^*\omega$ is a sum of commuting semipositive $(1,1)$-forms, so $\Omega^k=\sum_m\binom km\pi^*\omega_E^{k-m}\wedge p^*\omega^m$, each term being a non-negative measure on $\Gamma_f$. Divide by $k!$ and pull back by $(\mathrm{id},f)$.
\end{proof}

\begin{remark}[$V_k$ does not control the other terms uniformly]\label{rem:Vk}
Take $f_N\colon\Delta^2_*\to\C^2/\Lambda$ ($\Lambda$ a lattice, flat torus), $f_N(z)=[(Nz_1,z_2/N)]$. Then $f_N^*\omega=\frac i2(N^2dz_1\wedge d\bar z_1+N^{-2}dz_2\wedge d\bar z_2)$, so $(f_N^*\omega)^2=2\,dV$ is independent of $N$, whereas $\omega_E\wedge f_N^*\omega=(N^2+N^{-2})\,dV$. Each $f_N$ extends, so this shows only that no inequality $V_1\le C\,V_k$ with $C$ independent of $f$ can hold. Moreover $V_k\equiv0$ whenever the generic rank of $f$ is $<k$ (e.g.\ $\dim X<k$). Whether $V_k<\infty$ forces $V_1<\infty$ for a single map is Question~\ref{q:single}; Theorem~\ref{thm:pluri} answers it when $k=2$ for holomorphic maps and pluriclosed targets.
\end{remark}

\section{Bishop's extension theorem}\label{sec:bishop}

\begin{theorem}[Bishop \cite{Bis}; see \cite{Chi}]\label{thm:bishop}
Let $D\subset\C^p$ be open, $E\subsetneq D$ a closed analytic subset of dimension $<p$, $W\subset\C^m$ open, and $A$ a pure $p$-dimensional analytic subset of $(D\setminus E)\times W$ with locally finite $2p$-dimensional volume in $D\times W$, i.e.\ $\Hm^{2p}(A\cap L)<\infty$ for every compact $L\subset D\times W$. Then the closure $\ovl A$ of $A$ in $D\times W$ is a pure $p$-dimensional analytic subset of $D\times W$.
\end{theorem}

This and Lelong's theorem are the only external analytic inputs besides standard facts on proper mappings (Remmert, Grauert, Cartan--Remmert). The exceptional set lies in the base; the fibre $W$ is arbitrary.

\section{Proofs for maps}\label{sec:maps}

\begin{lemma}[Analytic closure]\label{lem:closure}
Let $Y$ be a complex manifold of dimension $n$, $Z\subset Y$ a closed analytic subset with $\dim Z<n$, $X$ a complex manifold with a Hermitian metric, and $\Gamma\subset(Y\setminus Z)\times X$ a pure $n$-dimensional analytic subset with $\Hm^{2n}(\Gamma\cap L)<\infty$ for every compact $L\subset Y\times X$. Then the closure $\ovl\Gamma$ of $\Gamma$ in $Y\times X$ is a pure $n$-dimensional analytic subset, and no irreducible component of $\ovl\Gamma$ lies in $Z\times X$.
\end{lemma}

\begin{proof}
Away from $Z\times X$, $\Gamma=\ovl\Gamma$. Let $(a,x)\in Z\times X$, $V\ni a$ a coordinate ball in $Y$ (relatively compact in a larger chart) and $W\ni x$ a coordinate ball in $X$ (likewise). The set $A=\Gamma\cap((V\setminus Z)\times W)$ is a pure $n$-dimensional analytic subset of $(V\setminus Z)\times W$ with finite $\Hm^{2n}$-measure on every compact subset of $V\times W$ (Euclidean and Hermitian metrics are comparable on compact sets). By Theorem~\ref{thm:bishop} ($p=n$, $E=Z\cap V$), $\ovl\Gamma\cap(V\times W)=\ovl A$ is analytic and pure $n$-dimensional. Analyticity is local, so $\ovl\Gamma$ is analytic. Finally, let $C$ be an irreducible component of $\ovl\Gamma$. The set $C\setminus\bigcup_{C'\ne C}C'$ is nonempty and open in $\ovl\Gamma$ (local finiteness of the family of components), hence meets the dense subset $\Gamma$; so $C\not\subset Z\times X$.
\end{proof}

\begin{lemma}[From the closed graph to a meromorphic map]\label{lem:tomero}
Let $Y$ be a complex manifold, $Z\subset Y$ a closed nowhere dense analytic subset, $f\colon Y\setminus Z\dashto X$ a meromorphic map with graph $\Gamma$ and $\ovl\Gamma$ the closure of $\Gamma$ in $Y\times X$. Assume that $\ovl\Gamma$ is analytic and $\pi\colon\ovl\Gamma\to Y$ is proper. Then $\ovl\Gamma$ is the graph of a meromorphic map $\bar f\colon Y\dashto X$ extending $f$.
\end{lemma}

\begin{proof}
\emph{Irreducibility.} If $\ovl\Gamma=A_1\cup A_2$ with $A_i$ closed analytic, then $\Gamma=(\Gamma\cap A_1)\cup(\Gamma\cap A_2)$, so $\Gamma\subset A_1$, say, by irreducibility of $\Gamma$; since $A_1$ is closed, $\ovl\Gamma\subset A_1$. Hence $\ovl\Gamma$ is irreducible, of dimension $n=\dim Y$. Moreover $\Gamma=\ovl\Gamma\cap((Y\setminus Z)\times X)$ because $\Gamma$ is closed in $(Y\setminus Z)\times X$.

\emph{Surjectivity.} By Remmert's proper mapping theorem $\pi(\ovl\Gamma)$ is a closed analytic subset of $Y$; it contains $\pi(\Gamma)=Y\setminus Z$ (the map $\pi|_\Gamma$ is proper with dense image), which is dense. So $\pi$ is onto.

\emph{The exceptional set.} Let $S\subset Y\setminus Z$ be the nowhere dense analytic set of Definition~\ref{def:mero} for $f$. By Grauert's direct image theorem $\pi_*\mathcal O_{\ovl\Gamma}$ is coherent. Since $\pi$ is surjective and $\ovl\Gamma$ is reduced, $\mathcal O_Y\to\pi_*\mathcal O_{\ovl\Gamma}$ is injective; its quotient $\mathcal Q$ is coherent and vanishes over $Y\setminus(Z\cup S)$, where $\pi$ is biholomorphic. Hence $S_1=\operatorname{supp}\mathcal Q\subset Z\cup S$ is a nowhere dense analytic subset of $Y$. Let $B=\{y\in\ovl\Gamma:\dim_y\pi^{-1}(\pi(y))\ge1\}$, an analytic subset of $\ovl\Gamma$ (Cartan--Remmert), proper since fibres over $Y\setminus(Z\cup S)$ are points; so $\dim B<n$ by irreducibility, and $S_2=\pi(B)$ is a nowhere dense analytic subset of $Y$. Put $S'=S_1\cup S_2$.

\emph{Biholomorphy.} Over $Y\setminus S'$ the proper map $\pi$ has finite fibres, hence is finite, and $\pi_*\mathcal O_{\ovl\Gamma}=\mathcal O_Y$; a finite map with this property is an isomorphism (the space $\ovl\Gamma$ is the analytic spectrum of $\pi_*\mathcal O_{\ovl\Gamma}$). So $\pi$ is biholomorphic over $Y\setminus S'$. Finally $\pi^{-1}(S')$ is a proper analytic subset of the irreducible space $\ovl\Gamma$ ($\pi$ onto, $S'$ nowhere dense), hence nowhere dense. Thus $\ovl\Gamma$ satisfies Definition~\ref{def:mero}. In particular the indeterminacy set $I=\{y\in Y:\dim\pi^{-1}(y)\ge1\}=S_2\subset S'$ has codimension $\ge2$: if $E$ were an irreducible component of $I$ of dimension $n-1$, then $\pi^{-1}(E)\to E$ would be onto with fibres of dimension $\ge1$, so $\dim\pi^{-1}(E)\ge n$, whereas $\pi^{-1}(E)\subset\pi^{-1}(S')$ is a nowhere dense analytic subset of the irreducible $n$-dimensional space $\ovl\Gamma$, hence of dimension $\le n-1$.
\end{proof}

\begin{lemma}[Local finiteness]\label{lem:localfin}
Let $X$ be a complex manifold with a Hermitian metric and $\Gamma\subset\Delta^k_*\times X$ a pure $k$-dimensional analytic subset such that $\pi\colon\Gamma\to\Delta^k_*$ is proper and $\Hm^{2k}(\Gamma\cap(U\times X))<\infty$ for a neighbourhood $U$ of $0$. Then $\Hm^{2k}(\Gamma\cap L)<\infty$ for every compact $L\subset\Delta^k\times X$.
\end{lemma}

\begin{proof}
Let $L\subset K\times X$ with $K\subset\Delta^k$ compact. Then $K\setminus U$ is a compact subset of $\Delta^k_*$, so $\pi^{-1}(K\setminus U)$ is compact in $\Gamma$ and has finite volume by Lelong's theorem. Hence $\Hm^{2k}(\Gamma\cap(K\times X))\le\Hm^{2k}(\Gamma\cap(U\times X))+\Hm^{2k}(\pi^{-1}(K\setminus U))<\infty$.
\end{proof}

\begin{proof}[Proof of Theorem~\ref{thm:maps}]
Let $\Gamma=\Gamma_f$. Shrinking $U$ we may assume $U\Subset\Delta^k$; by Lemma~\ref{lem:expansion}, \eqref{eq:fv} is equivalent to finiteness of $V_1,\dots,V_k$. By Lemma~\ref{lem:localfin} the hypotheses of Lemma~\ref{lem:closure} hold with $Y=\Delta^k$, $Z=\{0\}$, so $\ovl\Gamma$ is analytic and pure $k$-dimensional. Since $X$ is compact, $\pi\colon\ovl\Gamma\to\Delta^k$ is the restriction of the proper projection $\Delta^k\times X\to\Delta^k$ to a closed subset, hence proper. Lemma~\ref{lem:tomero} shows that $\ovl\Gamma$ is the graph of a meromorphic map extending $f$.
\end{proof}

\begin{proof}[Proof of Theorem~\ref{thm:bounded}]
Lemmas~\ref{lem:localfin} and~\ref{lem:closure} apply as above, so $\ovl\Gamma$ is analytic and $\ovl\Gamma\cap(U\times X)\subset U\times K$ ($K$ is closed). It remains to see that $\pi$ is proper. Let $L\subset\Delta^k$ be compact and $U'\Subset U$ a neighbourhood of $0$. Then $\pi^{-1}(L)\subset\pi^{-1}(L\setminus U')\cup\bigl((L\cap\ovl{U'})\times K\bigr)$, where $\pi^{-1}(L\setminus U')$ is compact (properness of $\pi$ over $\Delta^k_*$); so $\pi^{-1}(L)$ is a closed subset of a compact set. Lemma~\ref{lem:tomero} concludes. For correspondences see the proof of Theorem~\ref{thm:corr}.
\end{proof}

\begin{remark}[The case $k=1$]\label{rem:k1}
For $k=1$, \eqref{eq:fv} reduces to $\int f^*\omega<\infty$, the area counted with multiplicity: finite area implies extension. The hypothesis cannot be dropped: $f\colon\Delta_*\to E=\C/\Lambda$, $f(z)=[1/z]$, has an essential singularity (cluster set $E$) and infinite area, since $\int|z|^{-4}\,dA=\infty$.
\end{remark}

\section{Correspondences}\label{sec:corr}

\begin{lemma}[Connectedness of the complement]\label{lem:conn}
Let $M$ be a connected complex manifold and $S\subsetneq M$ a proper analytic subset. Then $M\setminus S$ is connected and dense.
\end{lemma}

\begin{proof}
The interior of $S$ is open, and closed in $M$: if $p$ lies in its closure, then $p\in S$, and near $p$ the set $S$ is the common zero set of finitely many holomorphic functions on a connected neighbourhood $N$ of $p$, all vanishing on a nonempty open subset of $N$, hence on $N$; so $p$ is interior. As $S\ne M$ and $M$ is connected, the interior is empty, so $M\setminus S$ is dense. Suppose $M\setminus S=U_1\sqcup U_2$ with $U_i$ open nonempty. The function equal to $1$ on $U_1$ and $0$ on $U_2$ is bounded holomorphic on $M\setminus S$, so by Riemann's extension theorem it extends holomorphically to $M$; it is continuous with values in $\{0,1\}$ on a dense set, hence everywhere, and constant by connectedness, a contradiction.
\end{proof}

\begin{proof}[Proof of Theorem~\ref{thm:corr}]
Let $\Gamma=\Gamma_F$ and $\ovl\Gamma$ its closure. By Lemmas~\ref{lem:localfin} and~\ref{lem:closure}, $\ovl\Gamma$ is analytic, pure $k$-dimensional, with no component in $\{0\}\times X$; every irreducible component of $\ovl\Gamma$ is the closure of a component of $\Gamma$. The map $\pi\colon\ovl\Gamma\to\Delta^k$ is proper as in the proofs of Theorems~\ref{thm:maps} and~\ref{thm:bounded}. We do not claim that $\pi$ is finite: $\pi^{-1}(0)$ may be positive dimensional (Example~\ref{ex:corr}).

Let $B=\{y\in\ovl\Gamma:\dim_y\pi^{-1}(\pi(y))\ge1\}$, let $R\subset\ovl\Gamma$ be the support of the relative differentials $\Omega^1_{\ovl\Gamma/\Delta^k}$ (the ramification locus), and put $S'=\pi(B\cup R\cup\Sing\ovl\Gamma)$. These sets are analytic. Each of $B$, $R$, $\Sing\ovl\Gamma$ has dimension $<k$ in every component $\ovl C$ of $\ovl\Gamma$: indeed $\ovl C$ is the closure of a component $C$ of $\Gamma$, $C\setminus\pi^{-1}(S_F)$ ($S_F$ the exceptional set of $F$) is nonempty open dense in $C$, and there $\pi$ is an unramified covering with smooth source. By Remmert's theorem $S'$ is a nowhere dense analytic subset of $\Delta^k$. Over $\Delta^k\setminus S'$ the map $\pi$ is proper, finite, unramified with smooth source, hence a finite covering, whose number of sheets $d'$ is constant since $\Delta^k\setminus S'$ is connected (Lemma~\ref{lem:conn}). At a point of the dense set $\Delta^k_*\setminus(S'\cup S_F)$ the fibre of $\ovl\Gamma$ equals the fibre of $\Gamma$, so $d'=d$.

Finally $\pi^{-1}(S')$ is nowhere dense in $\ovl\Gamma$: it is analytic, and if it contained a nonempty open subset of a component $\ovl C$, we could take it inside $\ovl C\setminus B$, where $\pi$ has finite fibres and hence the image has dimension $k$, contradicting $\dim S'<k$. Thus $\ovl\Gamma$ is a meromorphic correspondence of degree $d$.
\end{proof}

\section{Characterization and applications}\label{sec:charac}

\subsection{Necessity and general exceptional sets}

\begin{proposition}[Converse]\label{prop:converse}
If a meromorphic map (or correspondence) $f\colon\Delta^k_*\dashto X$, $X$ compact, extends meromorphically to $\Delta^k$, then \eqref{eq:fv} holds.
\end{proposition}

\begin{proof}
The extended graph $\ovl\Gamma$ is a pure $k$-dimensional analytic subset of $\Delta^k\times X$; for $U\Subset\Delta^k$, $\ovl\Gamma\cap(U\times X)$ is relatively compact, hence of finite $\Hm^{2k}$-measure by Lelong's theorem. As $\Gamma\subset\ovl\Gamma$, \eqref{eq:fv} holds.
\end{proof}

\begin{theorem}\label{thm:general}
Let $Y$ be a connected complex manifold of dimension $n$, $Z\subsetneq Y$ a proper closed analytic subset, $X$ compact, and $f\colon Y\setminus Z\dashto X$ meromorphic. If $\Vol(\Gamma_f\cap(K\times X))<\infty$ for every compact $K\subset Y$, then $f$ extends to a meromorphic map $Y\dashto X$; the converse holds as in Proposition~\ref{prop:converse}.
\end{theorem}

\begin{proof}
By Lemma~\ref{lem:conn}, $Y\setminus Z$ is connected and dense, and $\dim Z<n$. Every compact $L\subset Y\times X$ lies in some $K\times X$ with $K$ compact, so the hypothesis of Lemma~\ref{lem:closure} holds. Hence $\ovl\Gamma$ is analytic, and proper over $Y$ because $X$ is compact. Lemma~\ref{lem:tomero} concludes.
\end{proof}

\subsection{Energy criterion}
Let $Y$ be a complex manifold of dimension $n$ with Hermitian metric $\omega_Y$ and volume form $dV=\omega_Y^n/n!$, $X$ a complex manifold with metric $\omega$, and $f$ holomorphic on an open $W\subset Y$, with operator norm $\nrm{df}$.

\begin{lemma}[Density bound]\label{lem:density}
For measurable $W'\subset W$,
\[
\Vol\bigl(\Gamma(f|_{W'})\bigr)=\int_{W'}\prod_{i=1}^n(1+\lambda_i)\,dV\le\int_{W'}(1+\nrm{df}^2)^n\,dV,
\]
where $\lambda_i\le\nrm{df}^2$ are the eigenvalues of $f^*\omega$ with respect to $\omega_Y$.
\end{lemma}

\begin{proof}
Diagonalize $f^*\omega$ with respect to $\omega_Y$ at each point: $(\omega_Y+f^*\omega)^n/n!=\prod(1+\lambda_i)\,dV$; then use Wirtinger's equality. Finally $\lambda_i\le\nrm{df}^2$ because $f^*\omega(v,v)=|df\,v|_\omega^2\le\nrm{df}^2|v|^2_{\omega_Y}$.
\end{proof}

\begin{corollary}[Energy criterion]\label{cor:energy}
Let $Y$ be a connected complex manifold of dimension $n$, $Z\subset Y$ a proper closed analytic subset, $X$ compact, and $f\colon Y\setminus Z\dashto X$ meromorphic. If $\int_{K\setminus(Z\cup I_f)}\nrm{df}^{2n}dV<\infty$ for every compact $K\subset Y$, then $f$ extends. If $Z$ is a closed complex submanifold of codimension $c\ge1$ (e.g.\ a point, $c=n$), this holds whenever $\nrm{df}\le C\operatorname{dist}(\cdot,Z)^{-a}$ on $K\setminus(Z\cup I_f)$, with $C=C(K)$ and $0\le a<c/n$; for $Z=\{0\}\subset\Delta^n$ the condition is $a<1$.
\end{corollary}

\begin{proof}
Since $(1+t)^n\le2^{n-1}(1+t^n)$, Lemma~\ref{lem:density} bounds the volume of the graph over $K\setminus(Z\cup I_f)$ by $2^{n-1}\bigl(\Vol(K)+\int\nrm{df}^{2n}\bigr)<\infty$; the part over $I_f$ is a proper analytic subset of the irreducible $\Gamma_f$, hence has measure zero. Apply Theorem~\ref{thm:general}.

For the growth condition, since $Z$ is a submanifold, $\Vol(\{\operatorname{dist}(\cdot,Z)<\varepsilon\}\cap K)\le C'\varepsilon^{2c}$ for small $\varepsilon$. For $a>0$ and large $s$, $\{\operatorname{dist}^{-2an}>s\}=\{\operatorname{dist}<s^{-1/(2an)}\}$ has volume $\le C's^{-c/(an)}$, so $\int_K\operatorname{dist}^{-2an}dV=\int_0^\infty\Vol(\{\operatorname{dist}^{-2an}>s\}\cap K)\,ds<\infty$ as soon as $c/(an)>1$. For $a=0$ there is nothing to prove.
\end{proof}

\begin{remark}\label{rem:holo}
In the growth condition, $\nrm{df}\le C\operatorname{dist}(\cdot,Z)^{-a}$ on $(Y\setminus Z)\setminus I_f$ forces $I_f\cap(Y\setminus Z)=\emptyset$. Indeed, let $p\in I_f\cap(Y\setminus Z)$ and $B\ni p$ a coordinate ball with $\ovl B\cap Z=\emptyset$. On $B\setminus I_f$, $\nrm{df}\le M$. Since $I_f$ has real codimension $\ge4$, any two points of $B\setminus I_f$ can be joined by a piecewise linear path in $B\setminus I_f$ of length at most $|x-y|+\varepsilon$ for every $\varepsilon>0$; hence $f$ is Lipschitz on $B\setminus I_f$ for the Riemannian distance of $X$ and extends continuously to $p$. As $\Gamma_f$ is the closure of the graph over $B\setminus I_f$, the fibre $\pi^{-1}(p)$ is a single point, contradicting $\dim\pi^{-1}(p)\ge1$.
\end{remark}

\begin{remark}[Sharpness of the exponent for a point]\label{rem:sharp}
For $k=1$, $f(z)=[\log z/2\pi i]\in\C/(\mathbb Z+\tau\mathbb Z)$ is well defined (changing the branch of $\log$ adds $1$ after division by $2\pi i$), has $\nrm{df}=1/(2\pi|z|)$, and its cluster set at $0$ is the whole elliptic curve, so it does not extend: $a=1$ fails ($n=c=1$). For $Z=\{0\}\subset\Delta^n$, $n\ge2$, the Hopf map $f=p|_{\Delta^n_*}$ of Corollary~\ref{cor:hopf} (the case $g=\mathrm{id}$, which is of type (ii)) satisfies $\nrm{df}=c/\nrm{z}$ for the metric of Example~\ref{ex:hopf} and does not extend, so $a=1$ fails as well. The criterion is sufficient, not necessary: $[z_1:z_2]\colon\Delta^2_*\to\Pj^1$ has $\nrm{df}\sim\nrm{z}^{-1}$ and extends (with an indeterminacy point at $0$).
\end{remark}

\subsection{Hopf manifolds}

\begin{example}[Hopf surface]\label{ex:hopf}
Let $0<|\lambda|<1$, $H=(\C^2\setminus\{0\})/\langle z\mapsto\lambda z\rangle$, $p$ the quotient map and $f=p|_{\Delta^2_*}$. Every orbit enters every ball $B_\varepsilon$, so the cluster set of $f$ at $0$ is all of $H$. If $f$ extended with graph $\ovl\Gamma$, then $\ovl\Gamma\supset\{0\}\times H$, of dimension $2$, whereas $\pi^{-1}(0)$ is a proper analytic subset of the irreducible $2$-dimensional $\ovl\Gamma$, hence of dimension $\le1$: a contradiction. The form $\omega_0=i\sum dz_j\wedge d\bar z_j/\nrm{z}^2$ is invariant under $z\mapsto\lambda z$ and descends to a Hermitian metric $\omega_H$; then $f^*\omega_H^2=8\nrm{z}^{-4}dV$ and $\int_{B_r}f^*\omega_H^2=16\pi^2\int_0^r\rho^{-1}d\rho=\infty$ (logarithmic divergence). So $V_2=\infty$: finite volume is not automatic for compact targets, and extension is not automatic for Hopf surfaces.
\end{example}

\begin{corollary}[Hopf manifolds]\label{cor:hopf}
Let $n\ge2$, $\gamma$ a linear contraction of $\C^n$, $H=(\C^n\setminus\{0\})/\langle\gamma\rangle$ and $f\colon\Delta^n_*\to H$ holomorphic. Then $f=p\circ g$ with $g\colon\Delta^n\to\C^n$ holomorphic, and exactly one of the following holds:
\begin{enumerate}
\item[(i)] $g(0)\ne0$: $f$ extends holomorphically to a neighbourhood of $0$, with finite graph volume and energy near $0$;
\item[(ii)] $g(0)=0$: $f$ does not extend meromorphically, and its graph has infinite volume and infinite energy near $0$.
\end{enumerate}
\end{corollary}

\begin{proof}
$\Delta^n_*\simeq S^{2n-1}$ and $\C^n\setminus\{0\}$ are simply connected for $n\ge2$, so $p$ is the universal covering and $f$ lifts to $g\colon\Delta^n_*\to\C^n\setminus\{0\}$, which extends to $\Delta^n$ by Hartogs' theorem. If $g(0)\ne0$, then $p\circ g$ is holomorphic near $0$, and finiteness of volume and energy follows from compactness of the graph over a small closed ball. If $g(0)=0$, the zero is isolated (as $g\ne0$ on $\Delta^n_*$), so $g$ is finite at $0$ and hence open: for every $\varepsilon$, $g(B_\varepsilon)\supset W$ for a neighbourhood $W$ of $0$, and, since $g(z)\ne0$ for $z\ne0$ and $g(0)=0$, $f(B_\varepsilon\setminus\{0\})\supset p(W\setminus\{0\})=H$, because every $\gamma$-orbit enters $W$. As in Example~\ref{ex:hopf}, extension is impossible; Proposition~\ref{prop:converse} gives infinite graph volume, and infinite energy follows from Lemma~\ref{lem:density} and $(1+t)^n\le2^{n-1}(1+t^n)$.
\end{proof}

\subsection{Targets with Stein universal cover}

\begin{proposition}\label{prop:stein}
Let $k\ge2$ and $X$ compact with universal cover $\widetilde X$ equal to $\C^n$ or to a Stein manifold. Then every holomorphic $f\colon\Delta^k_*\to X$ extends holomorphically to $\Delta^k$ (no volume hypothesis).
\end{proposition}

\begin{proof}
$\Delta^k_*$ is simply connected for $k\ge2$, so $f$ lifts to $g\colon\Delta^k_*\to\widetilde X\subset\C^N$ (a proper holomorphic embedding in the Stein case, so $\widetilde X$ is a closed complex submanifold). Extend each coordinate by Hartogs' theorem to $G\colon\Delta^k\to\C^N$. Since $\widetilde X$ is closed and $\Delta^k_*$ is dense, $G(\Delta^k)\subset\widetilde X$, and $G$ is holomorphic as a map into the submanifold $\widetilde X$; then $p\circ G$ extends $f$.
\end{proof}

The hypothesis $k\ge2$ is necessary (elliptic curve, $k=1$). In particular extension is automatic for primary Kodaira and Inoue surfaces (universal cover $\C^2$, resp.\ $\mathbb H\times\C$).

\subsection{A correspondence with positive-dimensional fibre}

\begin{example}\label{ex:corr}
Let $k=2$, $X=\Pj^1$ and $\Gamma=\{(z,[x:y])\in\Delta^2\times\Pj^1:\ z_2x^2=z_1y^2\}$, an irreducible hypersurface (hence of locally finite volume). Over $z_1z_2\ne0$ the fibre has two points ($(x/y)^2=z_1/z_2$); it collapses to the double point $[0:1]$ if $z_1=0\ne z_2$ and to $[1:0]$ if $z_2=0\ne z_1$. So $\Gamma_F=\Gamma\cap(\Delta^2_*\times\Pj^1)$ is a correspondence of degree $2$ (namely $\sqrt{\cdot}\circ[z_1:z_2]$) and Theorem~\ref{thm:corr} applies; the fibre of $\ovl{\Gamma_F}=\Gamma$ over $0$ is all of $\Pj^1$.
\end{example}

\section{Pluriclosed targets in dimension two}\label{sec:pluri}

Throughout this section $k=2$, $X$ is a compact complex manifold with a Hermitian metric $\omega$ such that $\ddc\omega=0$, where $\dc=\frac i2(\bd-\partial)$ and $\ddc=i\partial\bd$ (a pluriclosed, or SKT, metric). If $\dim X=2$ this is the Gauduchon condition, and such a metric exists on every compact complex surface by \cite{Gau}. Let $f\colon\Delta^2_*\to X$ be holomorphic and $T:=f^*\omega$, a smooth positive $(1,1)$-form on $\Delta^2_*$ with $\ddc T=0$ (pull-back commutes with $d$ and $\dc$). Write
\[
T=\frac i2\sum_{\alpha,\beta=1}^2\tab{\alpha}{\beta}\,dz_\alpha\wedge d\bar z_\beta,\qquad(\tab{\alpha}{\beta})\ \text{Hermitian positive semidefinite},
\]
so that $T\wedge\omega_E=(t^{1\bar1}+t^{2\bar2})\,dV$ with $dV$ Lebesgue measure. On $\C$ we use $\Delta=4\partial\bd=\partial_x^2+\partial_y^2$. We write $\partial_j=\partial/\partial z_j$ and $\bd_j=\partial/\partial\bar z_j$.

\begin{lemma}[Pluriclosed identity]\label{lem:pluri}
$\ddc T=0$ is equivalent to
\begin{equation}\label{eq:pluri}
\partial_1\bd_1t^{2\bar2}+\partial_2\bd_2t^{1\bar1}=\partial_2\bd_1t^{1\bar2}+\partial_1\bd_2t^{2\bar1}.
\end{equation}
\end{lemma}

\begin{proof}
$\partial\bd T=\frac i2\sum\partial_\delta\bd_\gamma\tab{\alpha}{\beta}\,dz_\delta\wedge d\bar z_\gamma\wedge dz_\alpha\wedge d\bar z_\beta$. Only index choices with $\{\delta,\alpha\}=\{\gamma,\beta\}=\{1,2\}$ contribute, and reordering to $\Theta=dz_1\wedge d\bar z_1\wedge dz_2\wedge d\bar z_2$ gives the signs $(\delta,\alpha,\gamma,\beta)=(1,2,1,2):+$, $(1,2,2,1):-$, $(2,1,1,2):-$, $(2,1,2,1):+$. Hence
\[
\partial\bd T=\frac i2\bigl(\partial_1\bd_1t^{2\bar2}-\partial_1\bd_2t^{2\bar1}-\partial_2\bd_1t^{1\bar2}+\partial_2\bd_2t^{1\bar1}\bigr)\Theta,
\]
which vanishes if and only if \eqref{eq:pluri} holds.
\end{proof}

Fix $r_1\in(0,1)$ and define the slice area
\[
a(z_1):=\int_{\{z_1\}\times\Delta(r_1)}T=\int_{\Delta(r_1)}t^{2\bar2}(z_1,z_2)\,dA(z_2)\ \ge0,\qquad z_1\in\Delta(r_1)\setminus\{0\},
\]
(only the $(2,2)$ coefficient survives on the slice since $dz_1=d\bar z_1=0$ there).

\begin{lemma}[Slice potential]\label{lem:slice}
The function $a$ is $C^\infty$ on $\Delta(r_1)\setminus\{0\}$, and $\Delta a=\varphi$ there, where $\varphi$ extends to a function in $C^\infty(\Delta(r_1))$. Consequently there are $r_2\in(0,r_1)$ and $c_0\ge0$ with
\[
a(z_1)=c_0\log\frac1{|z_1|}+O(1)\quad(z_1\to0),\qquad\text{in particular}\quad\int_{\Delta(r_2)}a\,dV_{z_1}<\infty.
\]
\end{lemma}

\begin{proof}
For $z_1$ in a compact subset of $\Delta(r_1)\setminus\{0\}$, the closed slice $\{z_1\}\times\ovl{\Delta(r_1)}$ lies in $\Delta^2_*$, where $\tab{\alpha}{\beta}$ is smooth; hence $a\in C^\infty(\Delta(r_1)\setminus\{0\})$ and differentiation under the integral is legitimate:
\[
\Delta a=4\int_{\Delta(r_1)}\partial_1\bd_1t^{2\bar2}\,dA(z_2).
\]
Insert \eqref{eq:pluri}, written as $\partial_1\bd_1t^{2\bar2}=\partial_2(\bd_1t^{1\bar2})+\bd_2(\partial_1t^{2\bar1})-\partial_2(\bd_2t^{1\bar1})$. Each term on the right is a $\partial_2$- or $\bd_2$-derivative, and by Stokes in the fibre, with $dA=\frac i2dz_2\wedge d\bar z_2$,
\[
\int g_{,2}\,dA=\frac i2\oint g\,d\bar z_2,\qquad\int g_{,\bar2}\,dA=-\frac i2\oint g\,dz_2 .
\]
Applying this to $g=\bd_1t^{1\bar2}$, $g=\partial_1t^{2\bar1}$ and $g=\bd_2t^{1\bar1}$ yields
\begin{equation}\label{eq:phi}
\Delta a(z_1)=\varphi(z_1):=2i\oint_{|z_2|=r_1}\Bigl[(\bd_1t^{1\bar2})\,d\bar z_2-(\partial_1t^{2\bar1})\,dz_2-(\bd_2t^{1\bar1})\,d\bar z_2\Bigr].
\end{equation}
The integrand involves first derivatives of the smooth coefficients at points $(z_1,z_2)$ with $|z_2|=r_1>0$; these points lie in a compact subset of $\Delta^2_*$ for $z_1$ ranging in any compact subset of $\Delta(r_1)$, including $z_1=0$. So $\varphi\in C^\infty(\Delta(r_1))$.

Choose $r_2<r_1$ and $\psi\in C^\infty(\Delta(r_2))$, bounded, with $\Delta\psi=\varphi$ (Newtonian potential of a compactly supported smooth extension of $\varphi$). Then $h:=a-\psi$ is harmonic on $\Delta(r_2)\setminus\{0\}$ and $h\ge-\sup|\psi|$ since $a\ge0$. By B\^ocher's theorem \cite[Thm.~3.9]{ABR} applied to the positive harmonic function $h+\sup|\psi|$ on the punctured disc, $h=c_0\log\frac1{|z_1|}+u$ with $c_0\ge0$ and $u$ harmonic on $\Delta(r_2)$. Thus $a=h+\psi=c_0\log\frac1{|z_1|}+O(1)$, which is integrable.
\end{proof}

\begin{proposition}\label{prop:V1}
Under the assumptions of this section, there is $\rho>0$ with $\int_{\Delta(\rho)^2\setminus\{0\}}\omega_E\wedge f^*\omega=V_1(\Delta(\rho)^2)<\infty$.
\end{proposition}

\begin{proof}
By Lemma~\ref{lem:slice} and Tonelli, $\int_{\Delta(r_2)\times\Delta(r_1)}t^{2\bar2}\,dV=\int_{\Delta(r_2)}a\,dV_{z_1}<\infty$. The hypotheses are symmetric under exchanging $z_1\leftrightarrow z_2$ (equation \eqref{eq:pluri} is invariant under this relabelling), so also $\int_{\Delta(r_1)\times\Delta(r_2')}t^{1\bar1}\,dV<\infty$ for some $r_2'$. Take $\rho=\min(r_2,r_2')$ and use $T\wedge\omega_E=(t^{1\bar1}+t^{2\bar2})\,dV$ together with $t^{1\bar1},t^{2\bar2}\ge0$.
\end{proof}

\begin{proof}[Proof of Theorem~\ref{thm:pluri}]
Let $U=\Delta(\rho)^2$. By Lemma~\ref{lem:expansion} (with $k=2$, $I_f=\emptyset$), $\Vol(\Gamma_f\cap(U\times X))=\frac12V_0+V_1+\frac12V_2$, where $V_0<\infty$ and $V_1<\infty$ by Proposition~\ref{prop:V1}. Hence \eqref{eq:fv} holds if and only if $V_2(U)=\int_{U\setminus\{0\}}(f^*\omega)^2<\infty$. Theorem~\ref{thm:maps} gives extension if $V_2<\infty$, and Proposition~\ref{prop:converse} gives $V_2<\infty$ if $f$ extends.
\end{proof}

\begin{remark}\label{rem:metric}
For compact $X$ any two Hermitian metrics are comparable, and for positive $(1,1)$-forms $0\le\alpha\le C\beta$ implies $\alpha^2\le C^2\beta^2$ on a surface. Hence, when $\dim X=2$, finiteness of $\int(f^*\omega)^2=\int f^*(\omega^2)$ does not depend on the choice of $\omega$: the criterion says that $f$ extends if and only if the image volume of $f$, counted with multiplicity, is finite.
\end{remark}

\begin{example}[Hopf surface, revisited]\label{ex:hopf2}
Let $X$, $H$, $f$ be as in Example~\ref{ex:hopf} and $\omega=\frac i2\,\partial\bd\nrm{z}^2/\nrm{z}^2=\frac12\omega_0$. Writing $u=\log\nrm{z}^2$, one has $\omega_0=i\partial\bd u+i\,\partial u\wedge\bd u$. The first term is $\ddc$-closed. For the second, $\bd(\partial u\wedge\bd u)=\bd\partial u\wedge\bd u$ and then
\[
\partial(\bd\partial u\wedge\bd u)=\bd\partial u\wedge\partial\bd u=-(\partial\bd u)^2=0,
\] because $\partial\bd u$ is the pull-back of the Fubini--Study form of $\Pj^1$, which has rank one on $\C^2\setminus\{0\}$. Hence $\ddc\omega=0$. Here $t^{2\bar2}=(|z_1|^2+|z_2|^2)^{-1}$ and
\[
a(z_1)=\int_{|z_2|<r_1}\frac{dA(z_2)}{|z_1|^2+|z_2|^2}=\pi\log\frac{|z_1|^2+r_1^2}{|z_1|^2}=2\pi\log\frac1{|z_1|}+O(1),
\]
so $V_1<\infty$ while $V_2=\infty$ (Example~\ref{ex:hopf}); $f$ does not extend, in agreement with Theorem~\ref{thm:pluri}.
\end{example}

\begin{remark}[Scope of Theorem~\ref{thm:pluri}]\label{rem:scope}
The proof uses only that $T$ is a smooth, positive, $\ddc$-closed $(1,1)$-form on $\Delta^2_*$ and that $f$ is holomorphic on $\Delta^2_*$. It does not treat meromorphic maps whose indeterminacy points (a discrete subset of $\Delta^2_*$) accumulate at $0$: then the slice area may be singular along the countable set $\pi(I_f)$, and extending the identity $\Delta a=\varphi$ across it would require $\ddc f^*\omega=0$ as a current across $I_f$, which we do not prove. It also does not apply when $k\ge3$ or when $\omega$ is only Gauduchon in dimension $\ge3$ ($\ddc\omega^{n-1}=0$ does not imply $\ddc\omega=0$).
\end{remark}

\section{Further consequences}\label{sec:applications}

Throughout this section $X$ is a compact complex manifold with a Hermitian metric $\omega$, and $V_m$, $\Gamma_f$, $I_f$ are as in Section~\ref{sec:prelim}. We record consequences of Theorems~\ref{thm:maps}, \ref{thm:corr} and~\ref{thm:general} which go beyond the model case of a punctured polydisc.

\subsection{Graph volume over a general base}

Let $Y$ be a complex manifold of dimension $n$ with a Hermitian metric $\omega_Y$, and let $f\colon W\dashto X$ be a meromorphic map defined on an open set $W\subset Y$.

\begin{lemma}\label{lem:volume-general}
For every measurable $K\subset W$,
\[
\Vol\bigl(\Gamma_f\cap(K\times X)\bigr)=\sum_{m=0}^{n}\frac{1}{m!\,(n-m)!}\int_{K\setminus I_f}\omega_Y^{\,n-m}\wedge(f^*\omega)^m,
\]
where $\Vol$ is the $2n$-dimensional Hausdorff measure for the product metric $\omega_Y+\omega$. If $X$ is compact, finiteness of the left-hand side does not depend on the choice of $\omega$ (nor, over a relatively compact $K$, of $\omega_Y$).
\end{lemma}

\begin{proof}
The expansion is proved exactly as Lemma~\ref{lem:expansion}: $\Omega=\pi^*\omega_Y+p^*\omega$ is a sum of commuting semipositive $(1,1)$-forms, Wirtinger's equality gives $\Hm^{2n}=\int\Omega^n/n!$ on $\Gamma_f$, and the part of $\Gamma_f$ over $I_f$ is a proper analytic subset of the irreducible set $\Gamma_f$, hence has measure zero. For the last assertion, if $\omega'\le C\omega$ then $\Omega'\le C'\Omega$ and hence $\Omega'^{\,n}\le C'^{\,n}\Omega^n$, since powers of positive $(1,1)$-forms are monotone; two Hermitian metrics on a compact manifold are comparable.
\end{proof}

\subsection{Correspondences over a general base}

Let $Y$ be a connected complex manifold of dimension $k$ with a Hermitian metric $\omega_Y$, and $Z\subsetneq Y$ a closed analytic subset. By Lemma~\ref{lem:conn} the set $Y\setminus Z$ is connected, so Definition~\ref{def:corr} makes sense with $D=Y\setminus Z$. Volumes in $Y\times X$ are taken for the product metric $\omega_Y+\omega$; over compact subsets of $Y$, finiteness does not depend on $\omega_Y$ (Lemma~\ref{lem:volume-general}).

\begin{theorem}\label{thm:corr-general}
Let $X$ be compact and $F\colon Y\setminus Z\dashto X$ a meromorphic correspondence of degree $d$ such that $\Vol(\Gamma_F\cap(K\times X))<\infty$ for every compact $K\subset Y$. Then $F$ extends to a meromorphic correspondence $\widehat F\colon Y\dashto X$ of degree $d$.
\end{theorem}

\begin{proof}
Put $D=Y\setminus Z$, $\Gamma=\Gamma_F\subset D\times X$, let $\ovl\Gamma$ be its closure in $Y\times X$, $\pi$ the projection to $Y$, and $S_F\subset D$ the exceptional set of Definition~\ref{def:corr} for $F$.

\emph{Step 1: the closure.} Since $Y$ is connected and $Z\ne Y$, $\dim Z<k$. Every compact $L\subset Y\times X$ lies in some $K\times X$ with $K\subset Y$ compact, so $\Hm^{2k}(\Gamma\cap L)<\infty$ by hypothesis. By Lemma~\ref{lem:closure} (with $n=k$), $\ovl\Gamma$ is a pure $k$-dimensional analytic subset of $Y\times X$ with no irreducible component in $Z\times X$, and $\ovl\Gamma\cap(D\times X)=\Gamma$ because $\Gamma$ is closed in $D\times X$. Thus $\Gamma$ is open in $\ovl\Gamma$. Since $X$ is compact, $\pi\colon\ovl\Gamma\to Y$ is proper.

\emph{Step 2: the bad set is small.} Let $B=\{y\in\ovl\Gamma:\dim_y\pi^{-1}(\pi(y))\ge1\}$, let $R$ be the support of $\Omega^1_{\ovl\Gamma/Y}$ and $\Sing\ovl\Gamma$ the singular locus; these are analytic subsets of $\ovl\Gamma$. Put
\[
E=\bigl(\ovl\Gamma\cap\pi^{-1}(Z)\bigr)\cup\bigl(\Gamma\cap\pi^{-1}(S_F)\bigr).
\]
If $p\in\ovl\Gamma\setminus E$, then $\pi(p)\in D\setminus S_F$, where $\Gamma\to D\setminus S_F$ is an unramified covering; as $\Gamma$ is open in $\ovl\Gamma$, $\ovl\Gamma$ is smooth at $p$ and $\pi$ is a local biholomorphism there, so $p\notin B\cup R\cup\Sing\ovl\Gamma$. Hence $B\cup R\cup\Sing\ovl\Gamma\subset E$. The set $E$ is nowhere dense in $\ovl\Gamma$: the first part is a proper analytic subset of every component (Step~1), and the second is nowhere dense in $\Gamma$ by Definition~\ref{def:corr}, hence in $\ovl\Gamma$. An analytic subset $A\subset\ovl\Gamma$ with $\dim A=k$ would contain a nonempty open subset of $\ovl\Gamma$ (a $k$-dimensional irreducible component of $A$ is a component of the pure $k$-dimensional $\ovl\Gamma$), which is impossible inside $E$. Therefore $B$, $R$ and $\Sing\ovl\Gamma$ all have dimension $<k$.

\emph{Step 3: the covering and its degree.} Put $S'=\pi(B\cup R\cup\Sing\ovl\Gamma)$. By Remmert's proper mapping theorem $S'$ is analytic, and $\dim S'<k$, so $S'$ is a proper analytic subset of $Y$. By Lemma~\ref{lem:conn}, $Y\setminus S'$ is connected and dense. Over $Y\setminus S'$ the map $\pi$ is proper, has finite fibres (as $B$ is removed), is unramified and has smooth source; a proper local homeomorphism with finite fibres is a covering, so $\pi$ is a finite covering there, with a number $d'$ of sheets constant by connectedness. The set $Z\cup S_F\cup S'$ is nowhere dense in $Y$, so there is a point $y\in Y\setminus(Z\cup S_F\cup S')$; the fibre of $\ovl\Gamma$ over $y$ equals the fibre of $\Gamma$, which has $d$ points. Hence $d'=d$.

\emph{Step 4: $\pi^{-1}(S')$ is nowhere dense.} It is analytic in $\ovl\Gamma$. If it contained a nonempty open subset $O$ of a component of $\ovl\Gamma$, we could shrink $O$ to avoid $B$ and the other components; then $\pi|_O$ has finite fibres, so $\dim\pi(O)=k$, while $\pi(O)\subset S'$ has dimension $<k$, a contradiction.

\emph{Conclusion.} $\ovl\Gamma$ is a pure $k$-dimensional reduced analytic subset of $Y\times X$, $\pi\colon\ovl\Gamma\to Y$ is proper, $S'\subset Y$ is a nowhere dense analytic set with $\pi^{-1}(S')$ nowhere dense in $\ovl\Gamma$, and $\pi$ is an unramified $d$-sheeted covering over $Y\setminus S'$. By Definition~\ref{def:corr} (with $D=Y$, connected), $\ovl\Gamma$ is the graph of a meromorphic correspondence $\widehat F$ of degree $d$, and $\widehat F=F$ on $D$ since $\ovl\Gamma\cap(D\times X)=\Gamma$.
\end{proof}

\begin{corollary}[Multivalued inverses of branched coverings]\label{cor:branched}
Let $q\colon W\to Y\setminus Z$ be a proper finite holomorphic map of degree $d$ from a pure $k$-dimensional reduced complex space, $h\colon W\to X$ holomorphic, and assume that $(q,h)\colon W\to(Y\setminus Z)\times X$ is injective outside a nowhere dense analytic subset of $W$. Put $\Gamma=(q,h)(W)$. If $\Vol(\Gamma\cap(K\times X))<\infty$ for every compact $K\subset Y$, then the correspondence $F=h\circ q^{-1}$ extends to a meromorphic correspondence $Y\dashto X$ of degree $d$.
\end{corollary}

\begin{proof}
Put $D=Y\setminus Z$ (connected by Lemma~\ref{lem:conn}) and let $A\subset W$ be a nowhere dense analytic subset such that $(q,h)$ is injective on $W\setminus A$. The map $(q,h)$ is proper and finite because $q$ is, so by Remmert's theorem $\Gamma$, with its reduced structure, is a closed pure $k$-dimensional analytic subset of $D\times X$, and $\pi|_\Gamma\colon\Gamma\to D$ is proper and finite (the preimage of a compact $K\subset D$ is $(q,h)(q^{-1}(K))$). Let $S_0=\pi(\Sing\Gamma\cup R)$, where $R$ is the support of $\Omega^1_{\Gamma/D}$. Since $\pi|_\Gamma$ is finite, it has generic rank $k$ on each component of $\Gamma$, so $\dim R<k$ and $\dim\Sing\Gamma<k$; thus $S_0$ is a nowhere dense analytic subset of $D$ (Remmert) and $\pi^{-1}(S_0)\cap\Gamma$ is nowhere dense in $\Gamma$ (finite map). Over $D\setminus S_0$, $\pi|_\Gamma$ is a finite unramified covering with smooth source, with a constant number of sheets $d'$ since $D\setminus S_0$ is connected (Lemma~\ref{lem:conn}). The set $q(A)$ is analytic of dimension $<k$, so there is $z\in D\setminus(S_0\cup q(A))$ over which $q$ has exactly $d$ preimages; these lie in $W\setminus A$, and $(q,h)$ maps them injectively onto $\Gamma\cap\pi^{-1}(z)$. Hence $d'=d$. Thus $\Gamma$ is the graph of a meromorphic correspondence of degree $d$ in the sense of Definition~\ref{def:corr}, and Theorem~\ref{thm:corr-general} applies.
\end{proof}

\newpage
\subsection{Algebraic targets}

\begin{corollary}[Rationality criterion]\label{cor:rational}
Let $Y$ be a smooth connected projective variety, $Z\subsetneq Y$ a proper subvariety, $X$ a projective manifold, and $f\colon Y\setminus Z\dashto X$ a meromorphic map. Equip $Y$ and $X$ with Fubini--Study metrics. Then $f$ is the restriction of a rational map $Y\dashto X$ if and only if $\Vol(\Gamma_f)<\infty$. In that case the indeterminacy set of the extension is an algebraic subset of codimension $\ge2$.
\end{corollary}

\begin{proof}
If $f$ is rational with algebraic graph $\widehat\Gamma\subset Y\times X$, then $\Gamma_f\subset\widehat\Gamma$, and $\widehat\Gamma$ has finite volume by Lelong's theorem since it is compact; so $\Vol(\Gamma_f)<\infty$. Conversely, if the volume is finite, then $\Vol(\Gamma_f\cap(K\times X))<\infty$ for every compact $K\subset Y$, and Theorem~\ref{thm:general} shows that $f$ extends to a meromorphic map $\widehat f\colon Y\dashto X$, whose graph $\widehat\Gamma$ is a closed analytic subset of the projective variety $Y\times X$ and is therefore algebraic by Chow's theorem. The projection $\pi\colon\widehat\Gamma\to Y$ is biholomorphic over the complement of a proper analytic subset $S'\subset Y$, which is algebraic by Chow's theorem; hence $\pi$ is a birational morphism and $\widehat f$ is a rational map. Its indeterminacy set is $\pi(B)$ with $B=\{y\in\widehat\Gamma:\dim_y\pi^{-1}(\pi(y))\ge1\}$ analytic, hence algebraic by Chow's theorem, and it has codimension $\ge2$ by Definition~\ref{def:mero}.
\end{proof}

By Lemma~\ref{lem:volume-general} the condition can be written as finiteness of the mixed integrals $\int_{Y\setminus(Z\cup I_f)}\omega_Y^{\,n-m}\wedge(f^*\omega)^m$, $1\le m\le n=\dim Y$.

\subsection{Families}

Let $P$ be a complex manifold of pure dimension $\tau$ with a Hermitian metric (we give $P\times\Delta^k\times X$ the product metric) and let $G\subset P\times\Delta^k_*\times X$ be a pure $(\tau+k)$-dimensional analytic subset, proper over $P\times\Delta^k_*$, with fibres $G_t=\Gamma_{F_t}$ ($t\in P$) graphs of meromorphic correspondences $F_t\colon\Delta^k_*\dashto X$.

\begin{proposition}\label{prop:families}
Assume that for every compact $L\subset P$ there is a neighbourhood $U$ of $0$ in $\Delta^k$ with $\Hm^{2(\tau+k)}\bigl(G\cap(L\times U\times X)\bigr)<\infty$. Then:
\begin{enumerate}
\item[(a)] the closure $\widehat G$ of $G$ in $P\times\Delta^k\times X$ is a pure $(\tau+k)$-dimensional analytic set with no irreducible component contained in $P\times\{0\}\times X$, and $\ovl{G_t}\subset\widehat G\cap(\{t\}\times\Delta^k\times X)$ for every $t$;
\item[(b)] for almost every $t\in P$ (with respect to the Riemannian measure of $P$) the graph $G_t$ has finite volume near $\{0\}\times X$, hence $F_t$ extends to a meromorphic correspondence $\widehat F_t$ of the same degree.
\end{enumerate}
\end{proposition}

\begin{proof}
(a) Apply Lemma~\ref{lem:closure} with $Y=P\times\Delta^k$, $Z=P\times\{0\}$ and $\Gamma=G$; the hypothesis away from $Z$ follows from Lelong's theorem as in Lemma~\ref{lem:localfin}. The inclusion is clear since $\ovl{G_t}\subset\widehat G$ and lies in the slice.

(b) Let $\mathrm{pr}\colon G\to P$ be the projection, which is $1$-Lipschitz for the product metric. Eilenberg's inequality \cite[2.10.25]{Fed} (with the Hausdorff measures $\Hm^{2k}$ on the fibres and $\Hm^{2\tau}$ on $P$) gives
\[
\int^{*}\Hm^{2k}\bigl(G\cap(L\times U\times X)\cap\mathrm{pr}^{-1}(t)\bigr)\,dt\le C\,\Hm^{2(\tau+k)}\bigl(G\cap(L\times U\times X)\bigr)<\infty,
\]
so the integrand is finite for almost every $t\in L$. Exhausting $P$ by countably many compact sets gives the claim, and Theorem~\ref{thm:corr} applies to such $t$.
\end{proof}

\begin{remark}
The slice $\widehat G_t$ may a priori be larger than $\ovl{G_t}$: it can contain additional $k$-dimensional components lying in $\{0\}\times X$ when $\dim X\ge k$. We do not claim that $t\mapsto\widehat F_t$ is holomorphic in $t$; see Question~\ref{q:families}. Note also that it is not clear that a uniform bound on $\Vol(G_t)$ controls the total volume of $G$, which involves the variation of $G_t$ with $t$; this is why Proposition~\ref{prop:families} assumes finiteness of the latter.
\end{remark}

\begin{remark}[Moduli spaces and period maps]\label{rem:moduli}
Let $X$ be a compact complex \emph{manifold} containing a moduli space $M$ as a
dense Zariski-open subset, and let $f\colon\Delta^k_*\to X$ be holomorphic with
values in $M$. By Theorem~\ref{thm:maps}, Proposition~\ref{prop:converse} and
Lemma~\ref{lem:volume-general}, $f$ extends meromorphically across $0$ if and only if
its graph has finite $2k$-dimensional volume near $\{0\}\times X$ for some---equivalently, for any---%
smooth Hermitian metric on the compact manifold $X$. Two caveats delimit the
scope of this observation.

\emph{(i) Smoothness of $X$.} The hypothesis that $X$ be a manifold is
restrictive. Several natural compactifications fail it: the Deligne--Mumford
compactification $\overline{M}_g$ is a compact complex orbifold rather than a
manifold, and the Baily--Borel compactification of a locally symmetric variety is
in general singular. The statement applies only when a smooth compactification is
at hand (for instance a smooth toroidal compactification), and the admissible
metric is then an arbitrary smooth Hermitian metric on $X$, with no intrinsic
relation to the geometry of $M$.

\emph{(ii) The Weil--Petersson metric.} The Weil--Petersson metric
$\omega_{\mathrm{WP}}$ is defined only on $M$; it is incomplete and, near the
boundary, in general not comparable to a smooth Hermitian metric of $X$, so it does
not extend to a Hermitian metric on $X$. It is therefore not
admissible in Theorems~\ref{thm:maps}--\ref{thm:bounded}, and finiteness of
$\int (f^*\omega_{\mathrm{WP}})^k$ is not among the quantities our criterion
controls; whether, for an admissible $\omega$, finiteness of the single integral
$\int (f^*\omega)^k$ suffices is Question~\ref{q:single}.

For period maps the natural comparison is not the present volume criterion but
Borel's extension theorem~\cite{Bo}, which extends a holomorphic map from a product
of punctured discs and discs into an arithmetic quotient of a bounded symmetric
domain, holomorphically into the Baily--Borel compactification; its proof rests on
the distance-decreasing property for invariant metrics. The asymptotics of general
period maps near the boundary are governed by Schmid's nilpotent orbit
theorem~\cite{Sc}. Our criterion is of a
different, purely volumetric nature and recovers neither; we make no claim
relating it to this Hodge-theoretic extension theory.
\end{remark}

\section{Open questions}\label{sec:open}

\begin{question}[Single integral]\label{q:single}
Let $X$ be compact, $k\ge2$ and $f\colon\Delta^k_*\dashto X$ meromorphic with $\int(f^*\omega)^k<\infty$. Does $f$ extend? Theorem~\ref{thm:pluri} gives a positive answer for $k=2$, $f$ holomorphic and $X$ carrying a pluriclosed metric (e.g.\ any compact surface). The case of meromorphic $f$ with indeterminacy accumulating at $0$, and the case $k\ge3$, are open. A positive answer in general would imply extension of every meromorphic map of generic rank $<k$ (for which $V_k=0$).
\end{question}

\begin{question}[Disk-convex targets]\label{q:disk}
A complex space $X$ is \emph{disk-convex} \cite{Iv04} if for every compact $K\subset X$ there is a compact $\widehat K\subset X$ such that every holomorphic $h\colon\Delta\to X$, continuous on $\ovl\Delta$, with $h(\partial\Delta)\subset K$ satisfies $h(\Delta)\subset\widehat K$. Let $X$ be disk-convex and $f$ meromorphic with finite graph volume near $\{0\}\times X$. Is the cluster set of $f$ at $0$ relatively compact (so that Theorem~\ref{thm:bounded} applies)?
\end{question}

\begin{question}[Growth threshold]\label{q:growth}
Let $V(r)=\Vol\bigl(\Gamma_f\cap(\{r<\nrm{z}<1\}\times X)\bigr)$. Bounded $V$ implies extension (Theorem~\ref{thm:maps}), and $V(r)\sim c\log(1/r)$ for the Hopf example. Does $V(r)=o(\log(1/r))$ imply extension?
\end{question}

\begin{question}[Higher dimensions]\label{q:higher}
Let $\dim X=n\ge3$ and $f\colon\Delta^n_*\dashto X$ meromorphic with an isolated singularity. For which closedness conditions on $\omega$ (pluriclosed, non-K\"ahler, \dots) does finiteness of a single mixed integral $V_m$ control the lower ones?
\end{question}

\begin{question}[Families]\label{q:families}
In the setting of Proposition~\ref{prop:families}, suppose only that $\sup_{t\in L}\Vol(G_t\cap(U\times X))<\infty$ for each compact $L\subset P$. Is $\widehat G$ analytic, and does $\widehat G_t=\ovl{G_t}$ hold for every $t$?
\end{question}

\section*{Acknowledgements}
The author thanks S.~Ivashkovich for foundational work and stimulating
questions that motivated this research.

\end{document}